\documentclass[12pt,a4paper]{article}%
\usepackage{graphicx}
\usepackage{amsmath}%
\usepackage{amsfonts}%
\usepackage{amssymb}
\newtheorem{theorem}{Theorem}

\newtheorem{thm}{Theorem}[section]

\newtheorem{cla}[thm]{Claim}
\newtheorem{lem}[thm]{Lemma}
\newenvironment{proof}[1][Proof]{\textbf{#1.} }{\ \rule{0.5em}{0.5em}}

\begin{document}

\title{On the asymptotic value of the choice number of complete multi-partite graphs}
\author{Nurit Gazit\thanks{The School of Mathematical Sciences, Tel Aviv University,
Tel Aviv 69978, Israel. e-mail: perfect@post.tau.ac.il.}, Michael
Krivelevich\thanks{The School of Mathematical Sciences, Tel Aviv
University, Tel Aviv 69978, Israel. e-mail:
krivelev@post.tau.ac.il. Research supported in part by USA-Israel
BSF Grant 2002-133, and by grant 64/01 from the Israel Science
Foundation.}}
\maketitle

\newpage

\begin{abstract}
We calculate the asymptotic value of the choice number of complete
multi-partite graphs, given certain limitations on the relation
between the sizes of the different sides. In the bipartite case,
we prove that if $n_0\le n_1$ and $\log n_0 \gg \log\log n_1$, then
$ch(K_{n_0,n_1}) = (1+o(1))\frac{\log_2{n_1}}{ \log_2{x_0}}$,
where  $x_0$ is the unique root of the equation $x -
1 - x^{\frac{k-1}{k}} = 0$ in the interval $[1, \infty)$ and $k =
\frac{ \log_2{n_1}}{ \log_2{n_0}}$. In the multipartite case, we
prove that if $n_0 \leq n_1 ... \leq n_s$, and $n_0$ is not too
small compared to $n_s$, then $ch(K_{n_0,...,n_s}) = (1+o(1))
\frac{\log_2{n_s}}{ \log_2{x_0}}$. Here $x_0$ is
the unique root of the equation $s x - 1 -
\sum_{j=0}^{s-1}{x^{\frac{k_j-1}{k_j}}} = 0$ in the interval $[1,
\infty)$, and for every $0 \leq i \leq s - 1$, $k_i = \frac{
\log_2{n_s}}{ \log_2{n_i}}$.

\end{abstract}

\bigskip


\newpage

\section{Introduction}


The \emph{choice number} $ch(G)$ of a graph $G = (V,E)$ is the
minimum number $k$ such that for every assignment of a list $S(v)$
of at least $k$ colors to each vertex $v\in V$, there is a proper
vertex coloring of $G$ assigning to each vertex $v$ a color from
its list $S(v)$. The concept of choosability was introduced by
Vizing in 1976 [2] and independently by Erd\H{o}s, Rubin and
Taylor in 1979 [1]. It is also shown in [1] that the choice number
of the complete bipartite graph $K_{n,n}$ satisfies $ch(K_{n,n}) =
(1 + o(1)) \log_2{n}$. In this paper we calculate the asymptotic
value of the choice number of a general complete bipartite graph
$K_{n_0, n_1}$ and then expand the result to the case of a
complete multi-partite graph. We begin by proving (note that
throughout this paper all logs are binary):

\begin{theorem}
Let $2 \leq n_0 \leq n_1$ be integers, and let $n_0 =
{(\log{n_1})}^{\omega(1)}$. Denote $k = \frac{ \log{n_1}}{
\log{n_0}}$. Let $x_0$ be the unique root of the equation $x - 1 -
x^{\frac{k-1}{k}} = 0$ in the interval $[1, \infty)$. Then
$ch(K_{n_0,n_1}) = (1+o(1))\frac{\log{n_1}}{ \log{x_0}}$.
\end{theorem}

As usual, $\omega(1)$ stands for a function tending to infinity
arbitrarily slowly as its variable tends to infinity.

We will prove the theorem in two parts, showing first the
upper bound and then the lower bound. In the graph $K_{n_0,n_1}$
we label the group of $n_0$ vertices by $V_0$ and the group of $n_1$
vertices by $V_1$.

\section{The Upper Bound}


\begin{theorem}

Let $2 \leq n_0 \leq n_1$ be integers. Denote $k = \frac{
\log{n_1}}{ \log{n_0}}$. Let $x_0$ be the unique root of the
equation $x - 1 - x^{\frac{k-1}{k}} = 0$ in the interval $[1,
\infty)$. Then $ch(K_{n_0,n_1}) \leq \lceil \frac{\log{n_1}}{
\log{x_0}} \rceil + 1$.

\end{theorem}

\begin{proof}

\begin{lem}
If there exists a $p, 0 \leq p \leq 1$, s.t. $n_0 p^r + n_1 (1-p)^r
\leq 1$ then $ch(K_{n_0,n_1}) \leq r $.
\end{lem}

\begin{proof}
We show that given, for each vertex $v \in V(K_{n_0,n_1})$, a set of
colors $S(v)$ of size $r$, there is a proper vertex coloring of
the graph, assigning to each vertex $v$ a color from $S(v)$.

We partition the set of all available colors $S = \bigcup_{v \in
V} S(v)$ into two subsets $S_1$ and $S_0$ in the following manner:
each color $c \in S$ is chosen randomly and independently with
probability $p$ to be in $S_1$, and with probability $1-p$ to be
in $S_0$. We will show that with positive probability the sets $S_0$
and $S_1$ chosen satisfy the condition: each vertex $v \in V_0$
has a color $c \in S(v)$ s.t. $c \in S_0$, and each vertex $v \in
V_1$ has a color $c \in S(v)$ s.t. $c \in S_1$. Given such $S_0$
and $S_1$, we can color each vertex in $V_0$ with a color from
$S_0$, and each vertex in $V_1$ with a color from $S_1$, and since
$S_0 \cap S_1 = \O$, we get a proper coloring.

For each $v \in V_1$ the probability that a bad event occurs, i.e.
that all the colors in $S(v)$ are chosen to be in $S_0$, is
$(1-p)^r$. For each $v \in V_0$ the probability that a bad event
occurs, i.e. that all the colors in $S(v)$ are chosen to be in
$S_1$, is $p^r$. Therefore the expectation of the number of bad
events that occur is $n_0 p^r + n_1{(1-p)}^r \leq 1$. Since either
$p > 0$ or $1 - p > 0$, we can assume w.l.o.g. that $1 - p > 0$.
Then since, for example, the case in which all the colors in $S$
are chosen to be in $S_0$ happens with probability ${(1-p)} ^
{|S|} > 0$, and gives $n_1$ bad events, the case in which $0$
events occur also happens with positive probability (otherwise the
expectation would be greater than 1). Therefore we get the
desirable partition. \end{proof}

\begin{lem}\label{le22}
Given $r$ s.t. ${(\frac{1}{n_0})}^{\frac{1}{r-1}} +
{(\frac{1}{n_1})}^{\frac{1}{r-1}}\ge 1$, let $p = \frac {
{(\frac{1}{n_0})}^{\frac{1}{r-1}}}{{(\frac{1}{n_0})}^{\frac{1}{r-1}}
+ {(\frac{1}{n_1})}^{\frac{1}{r-1}} } $. Then $n_0 p^r +
n_1{(1-p)}^r \leq 1$.
\end{lem}

\begin{proof}
If $p = \frac {
{(\frac{1}{n_0})}^{\frac{1}{r-1}}}{{(\frac{1}{n_0})}^{\frac{1}{r-1}}
+ {(\frac{1}{n_1})}^{\frac{1}{r-1}} } $ then
${(\frac{p}{1-p})}^{r-1} = \frac{n_1}{n_0}$. Therefore
\begin{eqnarray*}
n_0 p^r + n_1{(1-p)}^r &=& n_0 p^r + n_1 {(\frac{n_0}{n_1})}
{{p}^{r-1}}{(1 - p)} =  n_0 {p}^{r-1}\\
&=& n_0  {\left(\frac {
{(\frac{1}{n_0})}^{\frac{1}{r-1}}}{{(\frac{1}{n_0})}^{\frac{1}{r-1}}
+ {(\frac{1}{n_1})}^{\frac{1}{r-1}} }\right)}^{r-1} =  {\left({\frac {1}
{{(\frac{1}{n_0})}^{\frac{1}{r-1}} +
{(\frac{1}{n_1})}^{\frac{1}{r-1}}}} \right)} ^{r-1} \\
&\leq& 1\ .
\end{eqnarray*}
\end{proof}

All that remains now is to choose $r=r(n_0,n_1)$ satisfying the
condition of Lemma \ref{le22}. Let $r=\lceil \frac{\log{n_1}}{\log{x_0}} \rceil + 1$.
Then $ r - 1 \geq \frac{\log{n_1}}{\log{x_0}}$, and hence
$x_0 \geq n_1^{\frac{1}{r-1}}$.
Since the function $f_k(x) = x - 1 - x^{\frac{k-1}{k}}$, where $k
\geq 1$, is a monotonely increasing function in the interval
${[1,\infty)}$, and since $f_k(x_0) = 0$, it follows that $
n_1^{\frac {1}{r-1}} \leq 1 +
{n_1^{{\frac{1}{r-1}}{\frac{k-1}{k}}}} = 1 +
{(\frac{n_1}{n_0})}^{\frac{1}{r-1}} $ as required. \end{proof}

\section{The Lower Bound}


\begin{theorem}
If $2 \leq n_0 \leq n_1$ are integers, and $n_0 =
{(\log{n_1})}^{\omega(1)}$, then $ch(K_{n_0,n_1}) \geq (1-o(1))
\frac{\log{n_1}}{ \log{x_0}}$, where  $x_0$ is
the unique root of the equation $x - 1 - x^{\frac{k-1}{k}} = 0$ in
the interval $[1, \infty)$ and $k = \frac{ \log{n_1}}{
\log{n_0}}$.
\end{theorem}

\textbf{Proof.}

A {\em cover} of a hypergraph $H$ is a subset $M$ of the vertices of the
hypergraph such that every hyperedge of $H$ contains at least one
vertex of $M$. A minimum cover is a cover which has the least
cardinality among all covers.

Let us generate the hypergraph $H_0$ created by the color lists of
the vertices in $V_0$, i.e. the hypergraph whose vertices are the
colors $\bigcup_{v \in V_0} S(v)$, and whose edges are the lists
$S(v)$ for each $v \in V_0$. In the same way, we generate the
hypergraph $H_1$ created by the color lists of the vertices in
$V_1$.

For any $r$, if we wish to prove $ch(K_{n_0,n_1}) > r$, it is
enough to show that there are parameters $t\geq r$ and $0\leq l\leq t$
s.t. it is possible to choose for each vertex in $K_{n_0,n_1}$
a list of $r$ colors from $\{1,2,...t\}$, and the lists chosen
satisfy:

\begin{enumerate}
\item The minimum cover of the hypergraph $H_0$ created by the
color lists of the vertices in $V_0$ (i.e. the minimum size of a
set $L$ of colors s.t. for every $v \in V_0$, $S(v)$ contains at
least one of the colors in $L$) is of cardinality at least $l$ .
\item The minimum cover of the hypergraph $H_1$ created by the
color lists of the vertices in $V_1$ is of cardinality at least
$t-l+1$.
\end{enumerate}

If these conditions are satisfied, then when these color lists are
assigned to the vertices of $K_{n_0,n_1}$, the graph cannot be properly
colored.  This is because at
least $l$ colors are needed to color one side, and at least
$t-l+1$ to color the other. Since there are only $t$ colors in
all, at least one color will be chosen by both sides -- i.e., at
least two vertices on opposite sides must be given the same color,
implying that a proper coloring is not possible. Therefore, the choice
number of the graph is greater than $r$.

\begin{lem}
If there exist parameters $t$ and $l$ such that $t\geq r, 0\leq l\leq t$
and
\begin{equation}\label{star1}
2^t e^{{- \frac{{(l)}_r}{{(t)}_r}} n_1} + 2^t
e^{-{\frac{{(t-l)}_r}{{(t)}_r}} n_0} \leq 1
\end{equation}
then $ch(K_{n_0,n_1}) > r $.
\end{lem}

\begin{proof}
It is easy to see that at least $l$ colors are required for a cover of
the hypergraph $H_0$ created by the color lists of the vertices in $V_0$
if and only if for each subset $C$ of size
$t-l+1$ of $\{1,2,...t\}$ there is at least one $v \in V_0$ for
which $S(v) \subset C$. In the same way, the minimum cover of the
hypergraph $H_1$ created by
the color lists of the vertices in $V_1$ is at least $t-l+1$
if and only if for each subset $C$ of size $l$ of
$\{1,2,...t\}$ there is at least one $v \in V_1$ for which $S(v)
\subset C$.

For each vertex $v$ in $K_{n_0,n_1}$, let $S(v)$ be a random
subset of cardinality $r$ of $\{1,2,...t\}$, chosen uniformly and
independently among all ${t \choose r}$ subsets of cardinality
$r$ of $\{1,2,...t\}$. We wish to find an $r$ that guarantees that
with positive probability:
\begin{enumerate}
\item For every subset $C$ of size $t-l+1$ there is a vertex $v
\in V_0$ s.t. $S(v) \subset C$, and \item For every subset $C$ of
size $l$ there is a vertex $v \in V_1$ s.t. $S(v) \subset C$.
\end{enumerate}

To simplify the calculations, we will change Condition 1 above to the
stronger condition that:
\begin{enumerate}
\item For every subset $C$ of size $t-l$ there is a vertex $v \in
V_0$ s.t. $S(v) \subset C$.
\end{enumerate}

For each fixed subset $C$ of cardinality $l$ of $\{1,2,...t\}$ and
each $v \in V_1$, the probability that $S(v) \nsubseteq C$ is $1 -
\frac{l\cdot...\cdot (l-r+1)}{t \cdot...\cdot (t-r+1)} = 1 -
\frac{{(l)}_r}{{(t)}_r}$. Since there are $n_1$ vertices in $V_1$
and ${t \choose l}$ subsets of cardinality $l$ of $\{1,...t\}$,
and since the color groups of the vertices were chosen
independently, the probability that there is a subset $C$ of size
$l$ that does not contain $S(v)$ for any $v \in V_1$ is at most
${t\choose l}{\left(1 - \frac{{(l)}_r}{{(t)}_r}\right)}^{n_1} <$
$2^t e^{{-\frac{{(l)}_r}{{(t)}_r}} n_1}$. In a similar fashion,
the probability that there is a subset $C$ of size $t-l$ that does
not contain $S(v)$ for any $v \in V_0$ is at most ${t \choose
t-l}{\left(1 -\frac{{(t-l)}_r}{{(t)}_r}\right)}^{n_0} < 2^t
e^{-{\frac{{(t-l)}_r}{{(t)}_r}} n_0}$.

We are looking for an $r$ that guarantees that the probability
that at least one of Conditions 1 and 2 does not hold is
smaller than 1. Therefore it is enough to show the sum of these
probabilities is smaller than 1, i.e., it is enough to show: $2^t
e^{{- \frac{{(l)}_r}{{(t)}_r}} n_1} + 2^t
e^{-{\frac{{(t-l)}_r}{{(t)}_r}} n_0} \leq 1$.
\end{proof}

Before proceeding to find $t$ and $l$ required in Lemma 3.1,
we derive bounds on $x_0$ that will be useful at later
stages of the proof.

\begin{lem}
$2\le x_0(k) < \max(k, e+2)$
\end{lem}

\begin{proof}
We begin by showing that if $k > e+1$, then $x_0(k) < k$. Since
$f_k(x)=x-1-x^{\frac{k-1}{k}}$ is monotonely increasing, we need
to show that $f_k(k) >0$, or $k - k^{\frac{k-1}{k}} - 1 > 0$, or
${(k-1)}^{\frac{1}{k-1}} > k^{\frac{1}{k}}$. But the function
$h(x) = x^{\frac{1}{x}}$ is monotonely decreasing for $x > e$.
So if $k > e+1$ then $k - 1 > e$ and therefore
${(k-1)}^{\frac{1}{k-1}}
> k^{\frac{1}{k}}$.

It can easily be seen that $x_0$ increases monotonely as a
function of $k$ (i.e. if $k_2 \geq k_1$, $x_0(k_2) \geq
x_0(k_1)$). Therefore if $k \leq e+2$, then $x_0(k) \leq x_0(e+2)
< e+2$.

To prove the lower bound on $x_0$, observe that
$f_k(2)=2-1-2^{\frac{k-1}{k}}=1-2^{\frac{k-1}{k}}\le 0$ for every
$k\ge 1$.
\end{proof}

\begin{lem}\label{le33}
Let $n_0 = {(\log{n_1})}^{\omega(1)}$. Define $r_0 =
\frac{\log{n_1}}{\log{x_0}}$, $u = \frac{4
\log{\log{n_1}}}{\log{n_0}} r_0$ and $r = r_0 - u$. Then $r = (1 -
o(1))r_0$, and for $t = {(\frac{n_1}{n_0})}^{\frac{1}{r}}r^2$ and
$l = t \frac{1}{{(\frac{n_1}{n_0})}^\frac{1}{r} + 1}$, $2^t e^{{-
\frac{{(l)}_r}{{(t)}_r}} n_1} + 2^t
e^{-{\frac{{(t-l)}_r}{{(t)}_r}} n_0} \leq 1$.
\end{lem}

\begin{proof}
If $n_0 = {(\log{n_1})}^{\omega(1)}$ then $\log{\log{n_1}} \ll
\log{n_0}$, and therefore $u = o(r_0)$, and $r = (1 - o(1))r_0$,
as required. From the fact that $r = (1 - o(1))r_0$, it also
follows that $r = \omega(1)$. This is because $x_0 < \max(k,
e+2)$, and therefore, if $k \leq e+2$ then $r_0 =
\frac{\log{n_1}}{\log{x_0}} > \frac{\log{n_1}}{\log{(e+2)}} =
\omega(1)$, and otherwise  $r_0 = \frac{\log{n_1}}{\log{x_0}} >
\frac{\log{n_1}}{\log{k}} =
\frac{\log{n_1}}{\log{\frac{\log{n_1}}{\log{n_0}}}} =
\frac{\log{n_1}}{\log{\log{n_1}} - \log{\log{n_0}}} \geq
\frac{\log{n_1}}{\log{\log{n_1}}} = \omega(1)$. Hence $r = (1 -
o(1))r_0 = \omega(1)$.

Let us denote $l_0 = l$ and $l_1 = t - l$. Then $t - l_i =
t
\frac{{(\frac{n_1}{n_i})}^\frac{1}{r}}{{(\frac{n_1}{n_0})}^\frac{1}{r}
+ 1}$, and $2^t e^{{- \frac{{(l)}_r}{{(t)}_r}} n_1} + 2^t
e^{-{\frac{{(t-l)}_r}{{(t)}_r}} n_0} = \sum_{i=0}^{1} 2^t e^{-
\frac{{(t - l_i)}_r}{{(t)}_r} n_i}$. In order for this sum to be
not greater than $1$, it is enough to show that $\frac{{(t -
l_i)}_r}{{(t)}_r} n_i \gg t$ for $i = 0,1$. We begin by estimating
$\frac{{(t - l_i)}_r}{{(t)}_r} n_i$.

\begin{cla}
$\frac{{(t - l_i)}_r}{{(t)}_r} n_i >
\frac{1}{2e^2}
{\frac{n_1}{{\left({(\frac{n_1}{n_0})}^\frac{1}{r} + 1\right)}^r}}$
for $i = 0,1$.
\end{cla}

\begin{proof}
$\frac{{(t - l_i)}_r}{{(t)}_r} > {(\frac{t - l_i - r}{t -
r})}^r = {(\frac{t - l_i}{t})}^r {(\frac{t(t - l_i -
r)}{(t-l_i)(t-r)})}^r = {(\frac{t - l_i}{t})}^r {(1 - \frac{l_i
r}{(t - l_i)(t - r)})}^r
> {(\frac{t-l_i}{t})}^r {(1 - \frac{2 l_i r}{(t - l_i)t})}^r$,
where the last inequality is a result of $r < \frac{t}{2}$.

Now since $\frac{l_0 r}{(t - l_0)t} = \frac{l r}{(t - l)t} =
\frac{t \frac{1}{{(\frac{n_1}{n_0})}^\frac{1}{r}+1}r}{t^2
\frac{{(\frac{n_1}{n_0})}^\frac{1}{r}}{{(\frac{n_1}{n_0})}^\frac{1}{r}
+ 1}} = \frac{r}{t {(\frac{n_1}{n_0})}^\frac{1}{r}} \leq \frac{r
{(\frac{n_1}{n_0})}^\frac{1}{r}}{t} = \frac{1}{r} = o(1)$, and
$\frac{l_1 r}{(t - l_1)t} = \frac{(t - l) r}{l t} = \frac{r
{(\frac{n_1}{n_0})}^\frac{1}{r}}{t} = \frac{1}{r} = o(1)$ we get
(recalling that $1-x\ge e^{-x}/2$ for $0\le x\le 1/2$) $\frac{{(t
- l_i)}_r}{{(t)}_r} > {(\frac{t - l_i}{t})}^r \frac{1}{2e^2}$.
Therefore $\frac{{(t - l_i)}_r}{{(t)}_r} n_i > {(\frac{t -
l_i}{t})}^r n_i \frac{1}{2 e^2} =
\left(\frac{(\frac{n_1}{n_i})^{\frac{1}{r}}}
{(\frac{n_1}{n_0})^{\frac{1}{r}}+ 1}\right)^r n_i \frac{1}{2 e^2}
= \frac{1}{2e^2}
{\frac{n_1}{{\left({(\frac{n_1}{n_0})}^\frac{1}{r} +
1\right)}^r}}\,. $
\end{proof}

Hence in order to prove that (\ref{star1}) holds it is now enough
to prove that ${\frac{n_1}{{\left({(\frac{n_1}{n_0})}^\frac{1}{r}
+ 1\right)}^r}} \gg t$.

\begin{cla}
${\frac{n_1}{{\left({(\frac{n_1}{n_0})}^\frac{1}{r} + 1\right)}^r}} \gg
t$\,.
\end{cla}

\textbf{Proof.} ${\frac{n_1}{\left((\frac{n_1}{n_0})^{\frac{1}{r}} +
1\right)^r} = {\left(\frac{n_1 ^
{\frac{1}{r}}}{{(\frac{n_1}{n_0})}^\frac{1}{r} + 1}\right)}^r =
\left[
\frac{n_1^{\frac{1}{r_0}}} {{(\frac{n_1}{n_0})}^{\frac{1}{r_0}} +1}
\frac{n_1^{\frac{1}{r}-\frac{1}{r_0}}}
{(\frac{n_1}{n_0})^{\frac{1}{r}} + 1}
{({(\frac{n_1}{n_0})}^\frac{1}{r_0} + 1)} \right]}^r$\,.

Since $\frac{n_1 ^ {\frac{1}{r_0}}
}{\left({(\frac{n_1}{n_0})}^\frac{1}{r_0} + 1\right)} = \frac{n_1
^ {\frac{\log{x_0}} {\log{n_1}}} } { {(\frac{n_1}{n_0})} ^
{\frac{\log{x_0}}{\log{n_1}}} +1 }
 =
\frac{x_0}{ \frac{x_0}{x_0 ^ {\frac{\log{n_0}}{\log{n_1}}}} + 1 }
= \frac{x_0}{x_0 ^ {\frac{k - 1}{k}} +1} = 1$, we get

${\frac{n_1}{{\left({(\frac{n_1}{n_0})}^\frac{1}{r} + 1\right)}^r}}  =
{\left(n_1^ {\frac{1}{r} - \frac{1}{r_0}}
\frac{{(\frac{n_1}{n_0})}^\frac{1}{r_0} + 1}
{{(\frac{n_1}{n_0})}^\frac{1}{r} + 1} \right)}^r >
{\left(n_1 ^ {\frac{1}{r}
- \frac{1}{r_0}} \frac{{(\frac{n_1}{n_0})}^\frac{1}{r_0}}
{{(\frac{n_1}{n_0})}^\frac{1}{r}}\right)}^r $, where the last
inequality follows from $r < r_0$. So
${\frac{n_1}{{({(\frac{n_1}{n_0})}^\frac{1}{r} + 1)}^r}} >
{\left(n_1 ^{\frac{1}{r} - \frac{1}{r_0}}
\frac{{(\frac{n_1}{n_0})}^\frac{1}{r_0}}
{{(\frac{n_1}{n_0})}^\frac{1}{r}}\right)}^r = n_0 ^ {({\frac{1}{r} -
\frac{1}{r_0}})r} = n_0^{1 - \frac{r}{r_0} } = n_0^{\frac{u}{r_0}}
= n_0 ^ {\frac{4 \log {\log{n_1}}}{\log{n_0}}} = \log^{4}{n_1}$.

Let us now estimate $t = {(\frac{n_1}{n_0})}^{\frac{1}{r}} r^2$.
Observe that
$r^2 < {r_0}^2 = {(\frac{\log{n_1}}{\log{x_0}})}^2 \leq
\log^{2}{n_1}$.
Also,
$$
{\left(\frac{n_1}{n_0}\right)}^{\frac{1}{r}} = 2 ^ {\frac{\log{n_1} -
\log{n_0}}{r}} = 2 ^ {\frac{\log{n_1} - \log{n_0}}{{\left(1 - \frac{4
\log{\log{n_1}}}{\log{n_0}}\right)} \frac{\log{n_1}}{\log{x_0}} }} = x_0
^ \frac{\frac{\log{n_1} - \log{n_0}}{\log{n_1}}} {1 - \frac{4
\log{\log{n_1}}}{\log{n_0}}} \leq x_0^{1+o(1)}\,,
$$
where the last inequality stems from the assumption that $n_0 =
{(\log{n_1})}^{\omega(1)}$. Since $x_0 = O(k)$,
${(\frac{n_1}{n_0})}^{\frac{1}{r}} \leq x_0^{1+o(1)} =
{(O(k))}^{1+o(1)} = O({( \log{n_1} )}^ {1 + o(1)})$.
Therefore $t = {(\frac{n_1}{n_0})}^{\frac{1}{r}} r^2 = O((\log{n_1})
^ {3 + o(1)}) \ll \log^{4}{n_1}$. \end{proof}
\newline
This also ends the proof of Lemma \ref{le33}, and therefore of the
lower bound and of Theorem 1.

\section{Generalization - Multi-Partite Graphs}


We wish to estimate the choice number of a general $(s+1)$-partite
graph $K_{n_0,n_1,...,n_s}$. In the graph $K_{n_0,n_1,...,n_s}$ we
label the group of $n_i$ vertices by $V_i$, for each $0 \leq i \leq
s$. Using a proof similar to that of the bipartite case, we will
prove:

\begin{theorem}
Let $s\ge 1$ be a fixed integer. Let
$2 \leq n_0 \leq n_1 ... \leq n_s$, and assume that $n_0 =
{(\log{n_s})}^\alpha$, where $\alpha \geq 2
\sqrt{\frac{\log{n_s}}{\log{\log{n_s}}}} $. For every $0 \leq i
\leq s - 1$ denote $k_i = \frac{ \log{n_s}}{ \log{n_i}}$. Let
$x_0$ be the unique root of the equation $sx -
1 - \sum_{j=0}^{s-1}{x^{\frac{k_j-1}{k_j}}} = 0$ in the interval
$[1, \infty)$. Then $ch(K_{n_0,...,n_s})=
(1+o(1))\frac{\log{n_s}}{\log{x_0}}$.
\end{theorem}

Again we divide the proof into two parts -- the upper bound and the
lower bound.

\section{The Upper Bound for Multi-Partite Graphs}


\begin{theorem}

Let $2 \leq n_0 \leq ... \leq n_s$ be integers, and let $0 <
\epsilon < 1$ be a constant. For every $0 \leq i \leq s - 1$
denote $k_i = \frac{ \log{n_s}}{ \log{n_i}}$. Let $x_0$ be the
unique root of the equation $(s + \epsilon) \cdot x - 1 -
\sum_{j=0}^{s-1}{x^{\frac{k_j-1}{k_j}}} = 0$ in the interval $[1,
\infty)$. Define $r = \lceil \frac{\log{n_s}}{ \log{x_0}} \rceil +
1$. Then $ch(K_{n_0,...,n_s}) \leq r$, for $n_s$ large enough.

\end{theorem}

\begin{proof}

\begin{lem}
If there exist ${p_0,...p_s}$ such that $0 \leq p_i \leq 1$ for every
$0\leq i \leq s$, ${\sum_{i=0}^{s}} p_i = 1$ and
${\sum_{i=0}^{s}} n_i (1 - p_i)^r \leq 1$, then
$ch(K_{n_0,n_1,...,n_s}) \leq r $.
\end{lem}

\begin{proof}
The proof is identical to that of the bipartite case (Lemma 2.1),
only this time we partition the set of all available colors into
$s+1$ sets, using the probabilities $p_i$. A bad event for a
vertex $v \in V_i$ is one in which all the colors in $S(v)$ are
chosen to be in color groups other than $S_i$, and it happens with
probability $(1 - p_i)^r$. \end{proof}

\begin{lem}
Given $r$ s.t. $\sum_{i=0}^{s}{n_i^{-{\frac{1}{r-1}}}}\geq
s^{\frac{r}{r-1}}$, let $p_i = 1-\frac{s
n_i^{-{\frac{1}{r-1}}}}{\sum_{j=0}^{s} {n_j^{-{\frac{1}{r-1}}}}}$
for $0 \leq i \leq s$. Then $0 \leq p_i \leq 1$ for each $0 \leq i
\leq s$, ${\sum_{i=0}^{s}} p_i = 1$, and ${\sum_{i=0}^{s}} n_i (1
- p_i)^r \leq 1$.
\end{lem}

\begin{proof}
In order for $p_i$ to be non-negative, we must demand that for
every $0 \leq i \leq s$, $\frac{s
n_i^{-{\frac{1}{r-1}}}}{\sum_{j=0}^{s} {n_j^{-{\frac{1}{r-1}}}}}
\leq 1$, or $s \leq \sum_{j=0}^{s}
{({\frac{n_i}{n_j}})^{\frac{1}{r-1}}}$. But if $s ^
{\frac{r}{r-1}} \leq \sum_{j=0}^{s} {n_j^{-{\frac{1}{r-1}}}}$, then
for every $0 \leq i \leq s$, $s < s ^ {\frac{r}{r-1}} \leq
\sum_{j=0}^{s} {n_j^{-{\frac{1}{r-1}}}} \leq \sum_{j=0}^{s}
{({\frac{n_i}{n_j}})^{\frac{1}{r-1}}}$. Also,
$$
{\sum_{i=0}^{s}} p_i = s + 1 - {\sum_{i=0}^{s}} (1 - p_i) = s + 1
- {\sum_{i=0}^{s}} \frac{s
(n_i^{-{\frac{1}{r-1}}})}{\sum_{j=0}^{s} {n_j^{-{\frac{1}{r-1}}}}}
 = s + 1 - s = 1\,.
$$
If  $1 - p_i = \frac{s n_i^{-{\frac{1}{r-1}}}}{\sum_{j=0}^{s}
{n_j^{-{\frac{1}{r-1}}}}}$ then $ {(\frac{1-p_i}{1-p_j})}^{r-1} =
\frac{n_j}{n_i} $. Therefore, for any $i$,
\begin{eqnarray*}
{\sum_{j=0}^{s}} n_j{(1 - p_j)}^r &=& n_i{(1-p_i)}^{r-1}\sum_{j=0}^{s}{(1 - p_j)} =
s \cdot n_i {(1 - p_i)}^{r-1} \\
&=&  s \cdot n_i \left({\frac{s
n_i^{-{\frac{1}{r-1}}}}{\sum_{j=0}^{s} {n_j^{-{\frac{1}{r-1}}}}}
}\right)^{r-1} = \left({\frac{s^{\frac{r}{r-1}}}{\sum_{j=0}^{s}
{n_j^{-{\frac{1}{r-1}}}}} }\right)^{r-1}\\
 &\leq& 1\,.
\end{eqnarray*}
\end{proof}

Let $r = \lceil \frac{\log{n_s}}{\log{x_0}} \rceil + 1$. Then $ r
- 1 \geq \frac{\log{n_s}}{\log{x_0}}$, and thus $x_0 \geq n_s
^{\frac{1}{r-1}}$.

Since the function $g_{k_0,...k_{s-1},\epsilon}(x) = (s +
\epsilon) \cdot x - 1 - \sum_{j=0}^{s-1}{x^{\frac{k_j-1}{k_j}}}$,
where $k_j \geq 1$ for each $j$, is a monotonely increasing
function in the interval ${[1,\infty)}$, and since
$g_{k_0,...k_{s-1},\epsilon}(x_0) = 0$, it follows that for $r$
large enough, or for $n_s$ large enough (see Lemma \ref{le62}
below, and the beginning of the proof of Lemma \ref{le33}), $ s ^
{\frac{r}{r-1}} n_s^{\frac {1}{r-1}} \leq (s + \epsilon)
n_s^{\frac {1}{r-1}} \leq 1 + \sum_{j=0}^{s-1}{n_s^{\frac {1}{r-1}
\frac{k_j-1}{k_j}}} = 1+ \sum_{i=0}^{s-1}
{{(\frac{n_s}{n_i})}^{\frac{1}{r-1}}}$ as required. \end{proof}

\section{The Lower Bound for Multi-Partite Graphs}


\begin{theorem}
Let $2 \leq n_0 ... \leq n_s$ be integers, and let $n_0 =
{(\log{n_s})}^\alpha$, where $\alpha \geq 2
\sqrt{\frac{\log{n_s}}{\log{\log{n_s}}}} $. For every $0 \leq i
\leq s - 1$ denote $k_i = \frac{ \log{n_s}}{ \log{n_i}}$. Let
$x_0$ be the unique root of the equation $s \cdot x - 1 -
\sum_{j=0}^{s-1}{x^{\frac{k_j-1}{k_j}}} = 0$ in the interval $[1,
\infty)$. Then
$ch(K_{n_0,...,n_s}) \geq (1-o(1))\frac{\log{n_s}}{\log{x_0}}$.
\end{theorem}

\textbf{Proof.}
Similarly to the bipartite case, in order to
prove $ch(K_{n_0,...,n_s}) > r$, it is enough to show that there
are a $t \geq r$ and a sequence of $0 \leq l_i \leq t$ for which
$\sum_{i=0}^{s} l_i = t$, s.t. it is possible to choose for each
vertex in $K_{n_0,...,n_s}$ a list of $r$ of colors from
$\{1,2,...t\}$, and the lists chosen satisfy the following s
conditions: For each $0 \leq i \leq s-1$ the minimum cover of the
hypergraph created by the color lists of the vertices in $V_i$ is
of cardinality at least $l_i$, and the additional condition: the
minimum cover of the
hypergraph created by the color lists of the vertices in $V_s$ is
of cardinality at least $l_s +1$.

As in the bipartite case, if these conditions are satisfied, then
by the pigeonhole principle at least 2 vertices in different
groups must be given the same color, so the choice number is
greater than $r$.

\begin{lem}\label{le61}
If there exist a parameter $ t \geq r$ and a sequence of $0 \leq l_i \leq t$
for which $\sum_{i=0}^{s} l_i = t$ and
\begin{equation}\label{star2}
\sum_{i=0}^{s} 2^t e^{-{\frac{{(t-l_i)}_r}{{(t)}_r}} n_i} \leq 1
\end{equation}
then $ch(K_{n_0,...,n_s}) > r $.
\end{lem}

\begin{proof}
Similar to the bipartite case. \end{proof}

As in the bipartite case, we calculate bounds on $x_0$ that will
help us later on.

\begin{lem}\label{le62}
$\frac{s+1}{s}\le x_0 < \max(k_0, e+2)$
\end{lem}

\begin{proof}
Since for every $0 \leq i \leq s$, $n_0 \leq n_i$, it follows that
$k_0 =
\frac{\log{n_s}}{\log{n_0}} \geq \frac{\log{n_s}}{\log{n_i}} =
k_i$.
Therefore, for a given $x$ in the range $[1, \infty)$, $x ^
{\frac{k_0-1}{k_0}} \geq x ^ {\frac{k_i-1}{k_i}}$ for all $i$, and
$f_{k_0,\ldots k_{s-1}}(x)=s x - 1 - \sum_{i=0}^{s-1}{x ^ {\frac{k_i-1}{k_i}}} \geq s x - 1
-s x ^ {\frac{k_0-1}{k_0}} = s (x - x ^ {\frac{k_0-1}{k_0}}) - 1
\geq x - x ^ {\frac{k_0-1}{k_0}} - 1$ (note all these functions
increase monotonely as functions of $x$). Therefore the root
$x_0$ in the range $[1, \infty)$ of the first equation $s x - 1 -
\sum_{i=0}^{s-1}{x ^ {\frac{k_i-1}{k_i}}} = 0$, which is our
equation, is not greater than the root $x_1$ of the equation
$x - x ^{\frac{k_0-1}{k_0}} - 1 = 0$.

But the last equation is $f_{k_0}(x)=0$, and we already know from
the bipartite case that its root is smaller than $\max(k_0, e+2)$.

To prove the lower bound observe that
$f_{k_0,\ldots,k_{s-1}}(\frac{s+1}{s})=s+1-1-\sum_{j=0}^{s-1}
\left(\frac{s+1}{s}\right)^{\frac{k_j-1}{k_j}}\le s-s=0$, and thus by
monotonicity $x_0\ge \frac{s+1}{s}$.
\end{proof}

\begin{lem}\label{le63}
Let $n_0 = {(\log{n_s})}^\alpha$, where $\alpha \geq 2
\sqrt{\frac{\log{n_s}}{\log{\log{n_s}}}} $. Define $r_0 =
\frac{\log{n_s}}{\log{x_0}}$, $u = \frac{4
\log{\log{n_s}}}{\log{n_0}} r_0$ and $r = r_0 - u$. Then $r = (1 -
o(1))r_0$, and for $t = {(\frac{1}{s}
\sum_{j=0}^{s}{(\frac{n_s}{n_j})}^{\frac{1}{r}} - 1)} r^2$ and $t
- l_i = t \frac{s
{(\frac{n_s}{n_i})}^\frac{1}{r}}{\sum_{j=0}^{s}{(\frac{n_s}{n_j})}^\frac{1}{r}}$,
one has: $0 \leq l_i \leq t$,$\sum_{i=0}^{s} l_i = t$, and
$\sum_{i=0}^{s} 2^t e^{-{\frac{{(t-l_i)}_r}{{(t)}_r}} n_i} \leq
1$, i.e., the assumptions of Lemma \ref{le61} are satisfied.
\end{lem}

\textbf{Proof.} Since $n_0 = {(\log{n_s})}^{\omega(1)}$, it
follows that $r = (1 - o(1))r_0$, as in the bipartite case. Also,
again as in the bipartite case, from $x_0 < \max(k_0, e+2)$ it
follows that $r_0 = \omega(1)$, and therefore $r = \omega(1)$.

We need to show that for every $i$, $0 \leq l_i \leq t$, or $0
\leq t - l_i \leq t$. Since $t - l_i$ is obviously non-negative,
we need to prove that $t - l_i \leq t$, or $\frac{s
{(\frac{n_s}{n_i})}^\frac{1}{r}}{\sum_{j=0}^{s}{(\frac{n_s}{n_j})}^\frac{1}{r}}
\leq 1$, or $s \leq
\sum_{j=0}^{s}{(\frac{n_i}{n_j})}^\frac{1}{r}$. Since $n_0\le n_i$ for
every $i$, it is enough to show $s \leq
\sum_{j=0}^{s}{(\frac{n_0}{n_j})}^\frac{1}{r}$.

Since $r_0 = \frac{\log{n_s}}{\log{x_0}}$, we have: $x_0 = n_s
^{\frac{1}{r_0}}$, and so $s {n_s}^{\frac {1}{r_0}} = 1 +
\sum_{j=0}^{s-1}{n_s^{\frac {1}{r_0} \frac{k_j-1}{k_j}}} =
\sum_{j=0}^{s} {{(\frac{n_s}{n_j})}^{\frac{1}{r_0}}}$, or $s =
\sum_{j=0}^{s} {{(\frac{1}{n_j})}^{\frac{1}{r_0}}}$.
But
$$
\sum_{j=0}^{s}{\left(\frac{n_0}{n_j}\right)}^\frac{1}{r} =
\sum_{j=0}^{s}{{\left(\frac{1}{n_j}\right)}^\frac{1}{r_0}
\frac{n_0^{\frac{1}{r}}} {n_j ^ {\frac{1}{r} - \frac{1}{r_0}}}}
\geq \frac{n_0^{\frac{1}{r}}} {n_s ^ {\frac{1}{r} -
\frac{1}{r_0}}}
\sum_{j=0}^{s}{{\left(\frac{1}{n_j}\right)}^\frac{1}{r_0} } = s
\frac{n_0^{\frac{1}{r}}} {n_s ^ {\frac{1}{r} - \frac{1}{r_0}}}\,,
$$
so it is enough to show $\frac{n_0^{\frac{1}{r}}} {n_s ^
{\frac{1}{r} - \frac{1}{r_0}}} \geq 1$.
But
${\frac{1}{r} - \frac{1}{r_0}} =
\frac{1}{r} \frac{u}{r_0}$, so $\frac{1}{n_s ^ {\frac{1}{r} -
\frac{1}{r_0}}} = 2^ {-{\frac{1}{r} \log{n_s} \frac{u}{r_0}}} = 2
^ {-{\frac{1}{r} \log{n_s}  \frac{4 \log{\log{n_s}}}{\log{n_0}} }}
= 2 ^ {-{\frac{1}{r} \log{n_s} \frac{4}{\alpha}}} $. Also
$n_0^{\frac{1}{r}} = {(\log{n_s})} ^ {\alpha \frac{1}{r}} = 2 ^
{\frac{1}{r} \alpha \log{\log{n_s}}}$. Therefore
$$
\frac{n_0^{\frac{1}{r}}} {n_s ^ {\frac{1}{r} - \frac{1}{r_0}}} =
{(2 ^ {\alpha \log{\log{n_s}} - \log{n_s} \frac{4}{\alpha} })} ^
\frac{1}{r} \geq 1\,,
$$
where the last inequality stems from the condition on $\alpha$.
Also,
$$
\sum_{i=0}^{s} l_i = (s+1)t - \sum_{i=0}^{s}{(t-l_i)} = (s+1)t -
\sum_{i=0}^{s} {t \frac{s
{(\frac{n_s}{n_i})}^\frac{1}{r}}{\sum_{j=0}^{s}{(\frac{n_s}{n_j})}^\frac{1}{r}}
} = st + t - st = t\,.
$$

All that is left for us to verify is that Condition (\ref{star2})
is fulfilled. The proof is,
again, similar to the bipartite case.

\begin{cla}
$\frac{{(t - l_i)}_r}{{(t)}_r} n_i > \frac{{s} ^ {r}
n_s}{{\left(\sum_{j=0}^{s}{(\frac{n_s}{n_j})}^\frac{1}{r}\right)}^r}
\frac{1}{2 e^2}$ for $0 \leq i \leq s$.
\end{cla}

\begin{proof}
We have:
$\frac{{(t - l_i)}_r}{{(t)}_r} > {(\frac{t - l_i - r}{t - r})}^r =
{(\frac{t - l_i}{t})}^r {(1 - \frac{l_i r}{(t - l_i)(t - r)})}^r
> {(\frac{t-l_i}{t})}^r {(1 - \frac{2 l_i r}{(t - l_i)t})}^r$
where the last inequality is a result of $r < \frac{t}{2}$. By
definition $t - l_i = t \frac{s
{(\frac{n_s}{n_i})}^\frac{1}{r}}{\sum_{j=0}^{s}{(\frac{n_s}{n_j})}^\frac{1}{r}}$,
so $l_i = \frac{
t{\left(\sum_{j=0}^{s}{(\frac{n_s}{n_j})}^\frac{1}{r} - s
{(\frac{n_s}{n_i})}^\frac{1}{r}\right)}}
{\sum_{j=0}^{s}{(\frac{n_s}{n_j})}^\frac{1}{r}}$, and
$\frac{l_i}{t - l_i} =
\frac{\sum_{j=0}^{s}{(\frac{n_s}{n_j})}^\frac{1}{r} - s
{(\frac{n_s}{n_i})}^\frac{1}{r}} {s
{(\frac{n_s}{n_i})}^\frac{1}{r}} = \frac{1}{s}
\sum_{j=0}^{s}{(\frac{n_i}{n_j})}^\frac{1}{r} - 1$.

Now since $\frac{l_i r}{(t - l_i)t} = {(\frac{1}{s}
\sum_{j=0}^{s}{(\frac{n_i}{n_j})}^\frac{1}{r} - 1)} \frac{r}{t}
\leq {(\frac{1}{s} \sum_{j=0}^{s}{(\frac{n_s}{n_j})}^\frac{1}{r} -
1)} \frac{r}{t} = \frac{1}{r} = o(1)$, we get $\frac{{(t -
l_i)}_r}{{(t)}_r} > {(\frac{t - l_i}{t})}^r \frac{1}{2e^2}$.

Hence $\frac{{(t - l_i)}_r}{{(t)}_r} n_i > {(\frac{t -
l_i}{t})}^r n_i \frac{1}{2 e^2} =
 {\left(\frac{s
{(\frac{n_s}{n_i})}^\frac{1}{r}}{\sum_{j=0}^{s}{(\frac{n_s}{n_j})}^\frac{1}{r}}\right)}^r
n_i \frac{1}{2 e^2} = {\frac{s^r
{n_s}}{{\left(\sum_{j=0}^{s}{(\frac{n_s}{n_j})}^\frac{1}{r}\right)}^r}}
\frac{1}{2 e^2}$. \end{proof}

Therefore in order to prove that (\ref{star2}) holds it is now
enough to prove that ${\frac{s^r
{n_s}}{{\left(\sum_{j=0}^{s}{(\frac{n_s}{n_j})}^\frac{1}{r}\right)}^r}}
\gg t$ (assuming $s$ is constant).

\begin{cla}
${\frac{s^r
{n_s}}{{\left(\sum_{j=0}^{s}{(\frac{n_s}{n_j})}^\frac{1}{r}\right)}^r}} \gg
t$\ .
\end{cla}

\begin{proof} We have:
$$
{\frac{s^ r
n_s}{{\left(\sum_{j=0}^{s}{(\frac{n_s}{n_j})}^\frac{1}{r}\right)}^r}}
= {\left(\frac{s n_s ^
{\frac{1}{r}}}{\sum_{j=0}^{s}{(\frac{n_s}{n_j})}^\frac{1}{r}}\right)}^r
= {\left[ \frac{s n_s ^ {\frac{1}{r_0}}
}{\sum_{j=0}^{s}{(\frac{n_s}{n_j})}^\frac{1}{r_0}} \frac{n_s ^
{\frac{1}{r} -
\frac{1}{r_0}}}{\sum_{j=0}^{s}{(\frac{n_s}{n_j})}^\frac{1}{r}}
{\sum_{j=0}^{s}{\left(\frac{n_s}{n_j}\right)}^\frac{1}{r_0}}\right
]}^r\,.
$$
Since $\frac{s n_s ^ {\frac{1}{r_0}}
}{\sum_{j=0}^{s}{(\frac{n_s}{n_j})}^\frac{1}{r_0}} = \frac{s n_s ^
{\frac{\log{x_0}} {\log{n_s}}} } {
\sum_{j=0}^{s-1}{(\frac{n_s}{n_j})}^\frac{\log{x_0}} {\log{n_s}}
+1 }
= \frac{s x_0}{\sum_{j=0}^{s-1} {x_0 ^ {\frac{k_j - 1}{k_j}}} +1}
= 1$, we get ${\frac{s^ r
n_s}{{\left(\sum_{j=0}^{s}{(\frac{n_s}{n_j})}^\frac{1}{r}\right)}^r}}
= {\left(n_s ^ {\frac{1}{r} - \frac{1}{r_0}}
\frac{\sum_{j=0}^{s}{(\frac{n_s}{n_j})}^\frac{1}{r_0}}
{\sum_{j=0}^{s}{(\frac{n_s}{n_j})}^\frac{1}{r}}\right)}^r =
{\left(\frac{\sum_{j=0}^{s}{(\frac{1}{n_j})}^\frac{1}{r_0}}
{\sum_{j=0}^{s}{(\frac{1}{n_j})}^\frac{1}{r}}\right)}^r $. Now,
$$
\frac{\sum_{j=0}^{s}{(\frac{1}{n_j})}^\frac{1}{r_0}}
{\sum_{j=0}^{s}{(\frac{1}{n_j})}^\frac{1}{r}} =
\frac{\sum_{j=0}^{s}{(\frac{1}{n_j})}^\frac{1}{r}
(\frac{1}{n_j})^{\frac{1}{r_0} - \frac{1}{r}}}
{\sum_{j=0}^{s}{(\frac{1}{n_j})}^\frac{1}{r}} =
\frac{\sum_{j=0}^{s}{(\frac{1}{n_j})}^\frac{1}{r} n_j^{\frac{1}{r}
- \frac{1}{r_0}}} {\sum_{j=0}^{s}{(\frac{1}{n_j})}^\frac{1}{r}}
\geq \frac{\sum_{j=0}^{s}{(\frac{1}{n_j})}^\frac{1}{r}
n_0^{\frac{1}{r} - \frac{1}{r_0}}}
{\sum_{j=0}^{s}{(\frac{1}{n_j})}^\frac{1}{r}}\,,
$$
where the last inequality is a result of $n_i \geq n_0$ for all $1
\leq i \leq s$ and of $r < r_0$. So
$\frac{\sum_{j=0}^{s}{(\frac{1}{n_j})}^\frac{1}{r_0}}
{\sum_{j=0}^{s}{(\frac{1}{n_j})}^\frac{1}{r}} \geq
n_0^{\frac{1}{r} - \frac{1}{r_0}}$, and ${\frac{s^ r
n_s}{{\left(\sum_{j=0}^{s}{(\frac{n_s}{n_j})}^\frac{1}{r}\right)}^r}}
\geq {\left(n_0^{\frac{1}{r} - \frac{1}{r_0}}\right)}^r = n_0^{1 -
\frac{r}{r_0}} = n_0 ^ {\frac{u}{r_0}} = n_0 ^ {\frac{4 \log
{\log{n_s}}}{\log{n_0}}} = \log^{4}{n_s}$.

Let us now estimate $t = {(\frac{1}{s}
\sum_{j=0}^{s}{(\frac{n_s}{n_j})}^{\frac{1}{r}} - 1)} r^2$.
First,
$r^2 < {r_0}^2 = {(\frac{\log{n_s}}{\log{x_0}})}^2 \leq
{(\frac{\log{n_s}}{\log{\frac{s+1}{s}}})}^2 = C \log^2{n_s}$ where
$C = C(s)$ is a constant.
Second,
$$
{\left(\frac{n_s}{n_0}\right)}^{\frac{1}{r}} = 2 ^
{\frac{\log{n_s} - \log{n_0}}{r}} = 2 ^ {\frac{\log{n_s} -
\log{n_0}}{{\left(1 - \frac{4\log{\log{n_s}}}{\log{n_0}}\right)}
\frac{\log{n_s}}{\log{x_0}} }} = x_0 ^ \frac{\frac{\log{n_s} -
\log{n_0}}{\log{n_s}}} {1 - \frac{4 \log{\log{n_s}}}{\log{n_0}}}
\leq x_0^{1+o(1)}\,,
$$
where the last inequality stems from the assumption that $n_0 =
{(\log{n_s})}^{\omega(1)}$. Since $x_0 = O(k_0)$, we get:
${(\frac{n_s}{n_0})}^{\frac{1}{r}} \leq x_0^{1+o(1)} =
{(O(k_0))}^{1+o(1)} = O({( \log{n_s} )}^ {1 + o(1)})$.

Therefore
\begin{eqnarray*}
t &=& {\left(\frac{1}{s}
\sum_{j=0}^{s}{\left(\frac{n_s}{n_j}\right)}^{\frac{1}{r}} -
1\right)} r^2 = {\left(\frac{1}{s}
\sum_{j=0}^{s-1}{\left(\frac{n_s}{n_j}\right)}^{\frac{1}{r}} -
\frac{s - 1}{s}\right)} r^2\\
 &\leq& {\left(\frac{1}{s}
\sum_{j=0}^{s-1}{\left(\frac{n_s}{n_j}\right)}^{\frac{1}{r}}\right)}
r^2 \leq \frac{1}{s} s
{\left(\frac{n_s}{n_0}\right)}^{\frac{1}{r}} r^2 =
{\left(\frac{n_s}{n_0}\right)}^{\frac{1}{r}} r^2 = O({(\log{n_s})}
^ {3 +
o(1)})\\
& \ll& \log^{4}{n_s}\,.
\end{eqnarray*}
 \end{proof}
\newline
This also ends the proof of Lemma \ref{le63}, and therefore of the
lower bound of the multi-partite case and of Theorem 4.

\newpage

\end{document}